\documentclass{amsart}

\usepackage{color}
\usepackage{amsthm}
\usepackage{amsmath, amssymb}
\usepackage{latexsym}
\usepackage{indentfirst}
\usepackage{graphicx}
\usepackage{placeins}
\usepackage{booktabs}
\usepackage{multirow}

\theoremstyle{plain}

\theoremstyle{definition}

\theoremstyle{remark}

\DeclareMathOperator{\Tr}{Tr}
\DeclareMathOperator{\tr}{Tr}

\DeclareMathOperator{\dist}{dist}
\DeclareMathOperator{\sn}{sn}
\DeclareMathOperator{\cn}{cn}
\DeclareMathOperator{\dn}{dn}

\newcommand{\secref}[1]{Section~\ref{#1}}

\newcommand{\etal}{\textit{et al}{}}
\newcommand{\ie}{\textit{i.e.}{}}

\newcommand{\pole}{\text{pole}}

\newcommand{\ud}{\,\mathrm{d}}
\newcommand{\RR}{\mathbb{R}}

\newcommand{\ZZ}{\mathbb{Z}}
\newcommand{\Or}{\mathcal{O}}

\newcommand{\bvec}[1]{\boldsymbol{#1}}

\newcommand{\eps}{\epsilon}
\newcommand{\kB}{k_\mathrm{B}}
\newcommand{\eV}{\mathrm{eV}}

\newcommand{\abs}[1]{\lvert#1\rvert}

\renewcommand{\Re}{\mathrm{Re}~}
\renewcommand{\Im}{\mathrm{Im}~}
\newcommand{\BB}{\mathcal{B}}
\newcommand{\CC}{\mathbb{C}}

\begin{document}

\title{Pole-based approximation of the Fermi-Dirac function}

\author{Lin Lin}
\address{Program in Applied and Computational
  Mathematics, Princeton University, Princeton, NJ 08544.
  Email: linlin@math.princeton.edu}

\author{Jianfeng Lu}
\address{Program in Applied and Computational
  Mathematics, Princeton University, Princeton, NJ 08544.
  Email: jianfeng@math.princeton.edu}

\author{Lexing Ying}
\address{Department of Mathematics and ICES, University of
  Texas at Austin, 1 University Station/C1200, Austin, TX 78712.
  Email: lexing@math.utexas.edu}

\author{Weinan E}
\address{Department of Mathematics and PACM,
  Princeton University, Princeton, NJ 08544.
  Email: weinan@math.princeton.edu}

\centerline{Dedicated to Professor Andy Majda on the occasion of his 60th 
birthday}

\begin{abstract}
  Two approaches for the efficient rational approximation of the
  Fermi-Dirac function are discussed: one uses the contour integral 
  representation and conformal mapping
  and the other is based on a version of the multipole
  representation of the Fermi-Dirac function that uses only simple poles.
  Both representations have logarithmic computational
  complexity. They are of great interest for electronic structure
  calculations.
\end{abstract}

\maketitle

\section{Introduction}

Given an effective one-particle Hamiltonian $\bvec{H}$, the inverse
temperature $\beta = 1/\kB T$ and the chemical potential $\mu$, the
finite temperature single-particle density matrix of the system is
given by the Fermi operator
\begin{equation}\label{eq:FDdist}
\bvec{\rho} = {2}{(1 + \exp(\beta(\bvec{H}-\mu)))^{-1}} 
= 1 - \tanh\Bigl(\frac{\beta}{2}(\bvec{H}-\mu)\Bigr),
\end{equation}
where $\tanh$ is the hyperbolic tangent function.  

In the last decade or so, the development of accurate and numerically
efficient representations of the Fermi operator has attracted a great deal of
attention in the quest for linear scaling electronic structure methods
based on effective one-electron Hamiltonians.  These approaches have
a numerical cost that scales linearly with $N$, the number of electrons,
and thus hold the promise of making electronic structure analysis of
large systems feasible. Achieving linear scaling for realistic
systems is very challenging.  Formulations based on the Fermi
operator are appealing since this operator gives directly the single
particle density matrix without the need for diagonalizing the Hamiltonian.

From a computational viewpoint, one main issue is that the right hand
side of \eqref{eq:FDdist} is an operator-valued function.  To evaluate
this function, we have to replace or approximate it by something which
can be computed directly without diagonalization.  Obvious candidates
are polynomial or rational approximations.  Such an approach was first
introduced by Baroni and Giannozzi \cite{BaroniGiannozzi1992} and
Goedecker and co-workers \cite{goedeckercolombo1994,
  goedeckerteter1995} (see also the review article
\cite{goedecker1999}). Several improvements have been made since then,
for example, in \cite{parrinelloarxiv, parrinello2008, parrinello2005,
  headgordon2004, headgordon2003, Ozaki2007}.  These are put broadly
under the umbrella of ``Fermi operator expansion'' (abbreviated as
FOE).

From the viewpoint of efficiency, a major concern is the cost for
representing the Fermi operator as a function of $\beta \Delta E$ (for
finite temperature) or $\Delta E/E_g$ (for gapped systems) where
$\beta$ is the inverse temperature, $\Delta E$ is the spectral width
of the discretized Hamiltonian matrix and $E_g$ is the spectrum gap of
the Hamiltonian around the chemical potential.  Consider a finite
temperature gapless system for example, the cost of the FOE proposed
by Goedecker \etal\ scales as $\beta \Delta E$. The fast polynomial
summation technique introduced by Head-Gordon \textit{et al} 
\cite{headgordon2004, headgordon2003} reduces the cost to
$(\beta\Delta E)^{1/2}$. The cost of the hybrid algorithm proposed by
Parrinello \etal\ in a recent preprint \cite{parrinelloarxiv} scales
as $(\beta\Delta E)^{1/3}$. The cost was brought down to logarithmic
scaling $\ln(\beta\Delta E)$ in \cite{LinLuCarE2009} using a multipole
representation of the Matsubara expansion of the Fermi-Dirac function.

The purpose of this article is to introduce two alternative rational
expansions of the Fermi-Dirac function that use only simple poles and
have computational cost that scales logarithmically. 
The first strategy is to use the
contour integral and conformal mapping idea proposed recently in
\cite{HaleHighamTrefethen2008}. This will be presented in the next
section. The other strategy is to borrow ideas from
\cite{YingBirosZorin2004} and use a version of multipole expansion
\cite{LinLuCarE2009} that only involves simple poles. This will be
discussed in \secref{sec:FMM}.  Numerical examples illustrating the
efficiency and accuracy of the representations are discussed in
\secref{sec:results}.

\section{Rational expansions based on contour integral}

Our first approach is an adaptation of the ideas proposed recently in
\cite{HaleHighamTrefethen2008} based on contour integral
representation and conformal mapping.  Let us first briefly recall the
main idea of \cite{HaleHighamTrefethen2008}. Consider a function $f$
that is analytic in $\CC\backslash (-\infty,0]$ and an operator
$\boldsymbol{A}$ with spectrum in $[m, M] \subset \RR^+$, one wants to
evaluate $f(\boldsymbol{A})$ using a rational expansion of $f$ by
discretizing the contour integral
\begin{equation}
  f(\boldsymbol{A}) = \frac{1}{2\pi i} \int_{\Gamma} f(z) 
  (z-\boldsymbol{A})^{-1} \ud z.
\end{equation}
The innovative technique in \cite{HaleHighamTrefethen2008} was to
construct a conformal map that maps the stripe $S = [-K,K]\times[0,K']$
to the upper half (denoted as $\Omega^+$) of the domain $\Omega =
\CC\backslash \bigl((-\infty,0]\cup[m,M]\bigr)$. This special map from
$t \in S$ to $z \in \Omega^+$ is given by
\begin{equation}
  z = \sqrt{mM}\Bigl(\frac{k^{-1}+u}{k^{-1}-u}\Bigr), \quad 
  u = \sn(t) = \sn(t|k),\quad
  k = \frac{\sqrt{M/m} - 1}{\sqrt{M/m}+1}.
\end{equation}
Here $\sn(t)$ is one of the Jacobi elliptic functions and the numbers
$K$ and $K'$ are complete elliptic integrals whose values are given by
the condition that the map is from $S$ to $\Omega^+$.

Applying the trapezoidal rule with $Q$ equally spaced points in 
$(-K+iK'/2, K+iK'/2)$, 
\begin{equation}
  t_j = -K + \frac{iK'}{2} + 2\frac{(j-\tfrac{1}{2}) K}{Q},
  \quad 1\leq j\leq Q,
\end{equation}
we get the quadrature rule (denote $z_j = z(t_j)$)
\begin{equation}
  f_Q(\boldsymbol{A}) = \frac{-4K\sqrt{mM}}{\pi Qk}\Im 
  \sum_{j=1}^Q \frac{f(z_j) (z_j - \boldsymbol{A})^{-1} \cn(t_j)\dn(t_j)}
  {(k^{-1}-\sn(t_j))^2}.
\end{equation}
Here $\cn$ and $\dn$ are the other two Jacobi elliptic functions in
standard notation and the factor
$\cn(t_j)\dn(t_j)(k^{-1}-\sn(t_j))^{-2}\sqrt{mM}/k$ comes from the
Jacobian of the function $z(t)$.

It is proved in \cite{HaleHighamTrefethen2008} that the convergence is
exponential in the number of quadrature points $Q$ and the exponent
deteriorates only
logarithmically as $M/m \to \infty$:
\begin{equation}\label{eq:accuracy}
\lVert f(\boldsymbol{A}) - f_Q(\boldsymbol{A})\rVert
=\Or(e^{-\pi^2 Q/(\log(M/m)+3)}).
\end{equation}

To adapt the idea to our setting with the Fermi-Dirac function or the
hyperbolic tangent function, we face with two differences: First,
the $\tanh$ function has singularities on the imaginary axis. Second,
the operator we are considering, $\beta(\boldsymbol{H}-\mu)$, has
spectrum on both the negative and positive axis.

\subsection{Gapped case}

We first consider the case when the Hamiltonian $\boldsymbol{H}$ has a
gap in its spectrum around the chemical potential $\mu$, such that
$\dist(\mu, \sigma(\boldsymbol{H})) = E_g>0$. Physically, this will be
the case when the system is an insulator.

Let us consider $f(z) = \tanh(\tfrac{\beta}{2}z^{1/2})$ acting on the
operator $\boldsymbol{A} = (\boldsymbol{H}-\mu)^2$.  
Now, $f(z)$ has singularities only on $(-\infty, 0]$ and the spectrum
of $\boldsymbol{A}$ is contained in $[E_g^2, E_M^2]$, where 
\begin{equation*}
E_M = \max_{E\in\sigma(\boldsymbol{H})} \lvert E - \mu \rvert.
\end{equation*}
We note that obviously $E_M \leq \Delta E$.  Hence we are back in the
same scenario as considered in \cite{HaleHighamTrefethen2008}
except that we need to take care of different branches of the square
root function when we apply the quadrature rule.

More specifically, we construct the contour and quadrature points
$z_j$ in the $z$-plane using parameters $m = E_g^2$ and $M = E_M^2$.
Denote $g(\xi) = \tanh(\beta \xi/2)$, $\xi_j^{\pm} = \pm z_j^{1/2}$,
and $\boldsymbol{B} = \boldsymbol{H} - \mu$. The quadrature rule is
then given by
\begin{multline}
  g_Q(\boldsymbol{B}) = \frac{-2K\sqrt{mM}}{\pi Qk}\Im
  \Biggl(\sum_{j=1}^Q \frac{g(\xi_j^+) (\xi_j^+ - \boldsymbol{B})^{-1}
    \cn(t_j)\dn(t_j)}
  {\xi_j^+(k^{-1}-\sn(t_j))^2} \\
  + \sum_{j=1}^Q \frac{g(\xi_j^-) (\xi_j^- - \boldsymbol{B})^{-1}
    \cn(t_j)\dn(t_j)} {\xi_j^-(k^{-1}-\sn(t_j))^2}
  \Biggr),
\end{multline}
where the factors $\xi_j^{\pm}$ in the denominator come from the
Jacobian of the map from $z$ to $\xi$. The number of poles to be
inverted is $N_{\pole} = 2Q$. After applying \eqref{eq:accuracy}, we
have a similar error estimate for $g(\boldsymbol{B})$
\begin{equation}\label{eq:gapaccuracy}
  \lVert g(\boldsymbol{B}) - g_Q(\boldsymbol{B})\rVert
  =\Or(e^{-\pi^2 Q/(2\log(E_M/E_g)+3)}).
\end{equation}

\begin{figure}[ht]
  \begin{center}
    \includegraphics[width=8cm]{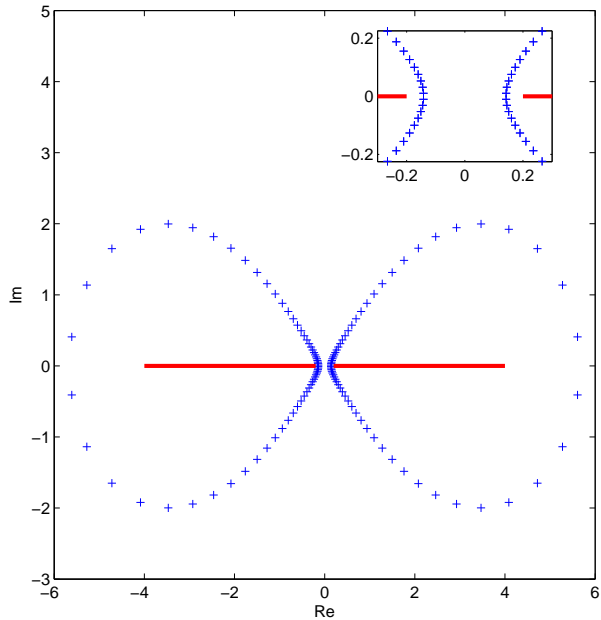}
  \end{center}
  \caption{A typical configuration of the poles on a two-loop contour.
    $Q = 30$, $E_g = 0.2$, $E_M = 4$ and $\beta = 1000$. The red line
    indicates the spectrum. The inset shows the poles close to the
    origin. The x-axis is $E-\mu$ with $E$ the eigenvalue of
    $\boldsymbol{H}$.  The poles with negative imaginary parts are not
    explicitly calculated.  }
  \label{fig:twocircle}
\end{figure}

In Fig.~\ref{fig:twocircle}, a typical configuration of the quadrature
points is shown. The x-axis is taken to be $E-\mu$. We see that in
this case the contour consists of two loops, one around the spectrum below
the chemical potential and the other around the spectrum above the
chemical potential.

Note further that as the temperature goes to zero, the
Fermi-Dirac function converges to the step function:
\begin{equation}
\eta(\xi) = 
\begin{cases}
2, & \xi\leq 0, \\
0, & \xi > 0.
\end{cases}
\end{equation}
Therefore, the contribution of the quadrature points $\xi_j^+$ on the
right half plane ($\Re\xi_j^+>0$) is negligible when $\beta$ is large.
In particular, for the case of zero temperature, one may choose only
the quadrature points on the left half plane. The quadrature formula
we obtain then becomes
\begin{equation}
\label{eqn:matrixsignquad}
\eta_Q(\boldsymbol{B}) = \frac{-4K\sqrt{mM}}{\pi Qk}\Im\biggl(
\sum_{j=1}^Q \frac{(\xi_j^- - \boldsymbol{B})^{-1}\cn(t_j)\dn(t_j)}{
\xi_j^{-}(k^{-1} - \sn(t_j))^2}\biggr).
\end{equation}
The number of poles to be inverted is then $N_{\pole}=Q$.

\begin{figure}[ht]
  \begin{center}
    \includegraphics[width=8cm]{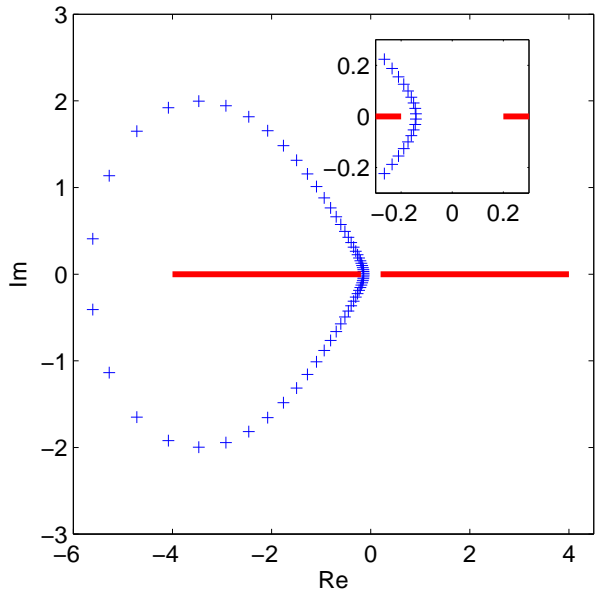}
  \end{center}
  \caption{A typical configuration of the poles for zero temperature
    ($\beta =\infty$). $Q = 30$, $E_g=0.2$ and $E_M = 4$. The red line
    indicates the spectrum. The inset zooms into the poles that is
    close to the origin.  The x-axis is $E-\mu$ with $E$ the
    eigenvalue of $\bvec{H}$. The poles with negative imaginary parts
    are not explicitly calculated.  }
  \label{fig:matrixsigncontour}
\end{figure}

We show in Fig.~\ref{fig:matrixsigncontour} a typical configuration of
the set of quadrature points. Only one loop is required compared with
Fig.~\ref{fig:twocircle}.

\subsection{Gapless case}

The more challenging case is when the spectrum of $\boldsymbol{H}$
does not have a gap, \textit{i.e.}, $E_g = 0$. Physically, this
corresponds to the case of metallic systems. In this case, the
construction discussed in the last subsection does not work.

\begin{figure}[ht]
  \begin{center}
    \includegraphics[width=8cm]{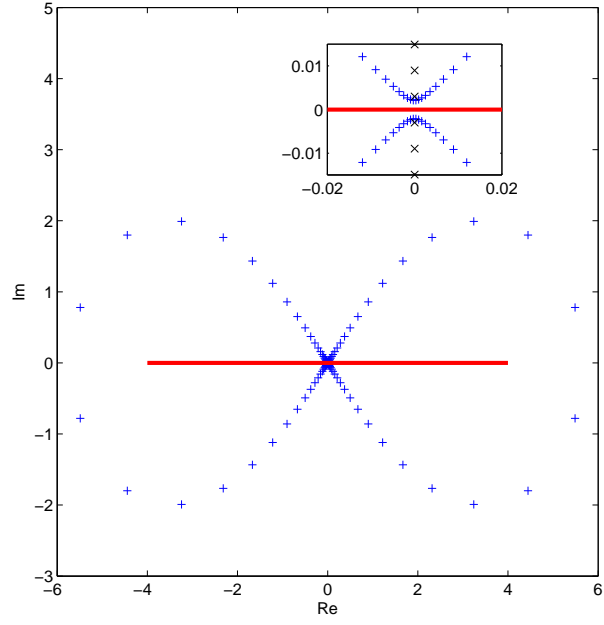}
  \end{center}
  \caption{A typical configuration of the poles on a dumbbell-shaped
    contour. $Q = 30$, $E_g = 0$, $E_M = 4$ and $\beta = 1000$.  The inset zooms
    into the part close to the origin. The red line indicates the
    spectrum. The black crosses indicate the positions of the poles of
    $\tanh$ function on the imaginary axis. The poles with negative
    imaginary parts are not explicitly calculated.  }
  \label{fig:dumbbell}
\end{figure}

To overcome this problem, we note that the hyperbolic tangent function
$\tanh(\tfrac{\beta}{2}z)$ is analytic except at poles
$(2l-1)\pi/\beta i,\ l\in\ZZ$ on the imaginary axis. Therefore, we
could construct a contour around the whole spectrum of
$\boldsymbol{H}$ which passes through the imaginary axis on the upper
half plane between the origin and $\pi/\beta i$ and also on the lower
half plane between the origin and $-\pi/\beta i$. Thus, we will have a
dumbbell shaped contour as shown in Fig.~\ref{fig:dumbbell}.

To be more specific, let us first construct the contour and quadrature
points $z_j$ in the $z$-plane as in the last subsection using
parameters $m=\pi^2/\beta^2$ and $M = E_M^2 + \pi^2/\beta^2$. Denote
$\xi_j^{\pm} = \pm (z_j - \pi^2/\beta^2)^{1/2}$, $g = \tanh(\beta
\xi/2)$ and $\boldsymbol{B} = \boldsymbol{H} - \mu$. The quadrature
rule takes the following form
\begin{multline}
  \label{eqn:gapless}
  g_Q(\boldsymbol{B}) = \frac{-2K\sqrt{mM}}{\pi Qk}\Im
  \Biggl(\sum_{j=1}^Q \frac{g(\xi_j^+) (\xi_j^+ - \boldsymbol{B})^{-1}
    \cn(t_j)\dn(t_j)}
  {\xi_j^+(k^{-1}-\sn(t_j))^2} \\
  + \sum_{j=1}^Q \frac{g(\xi_j^-) (\xi_j^- - \boldsymbol{B})^{-1}
    \cn(t_j)\dn(t_j)} {\xi_j^-(k^{-1}-\sn(t_j))^2}
  \Biggr).
\end{multline}
When apply the quadrature formula, the number of poles to be inverted is
$N_{\pole} = 2Q$. 
Fig.~\ref{fig:dumbbell} shows a typical configuration of quadrature
points for $Q=30$. The map $\xi(z) = (z-\pi^2/\beta^2)^{1/2}$ maps the
circle in the $z$-plane to a dumbbell-shaped contour (put two branches
together).

\begin{figure}[ht]
  \begin{center}
    \includegraphics[width=\textwidth]{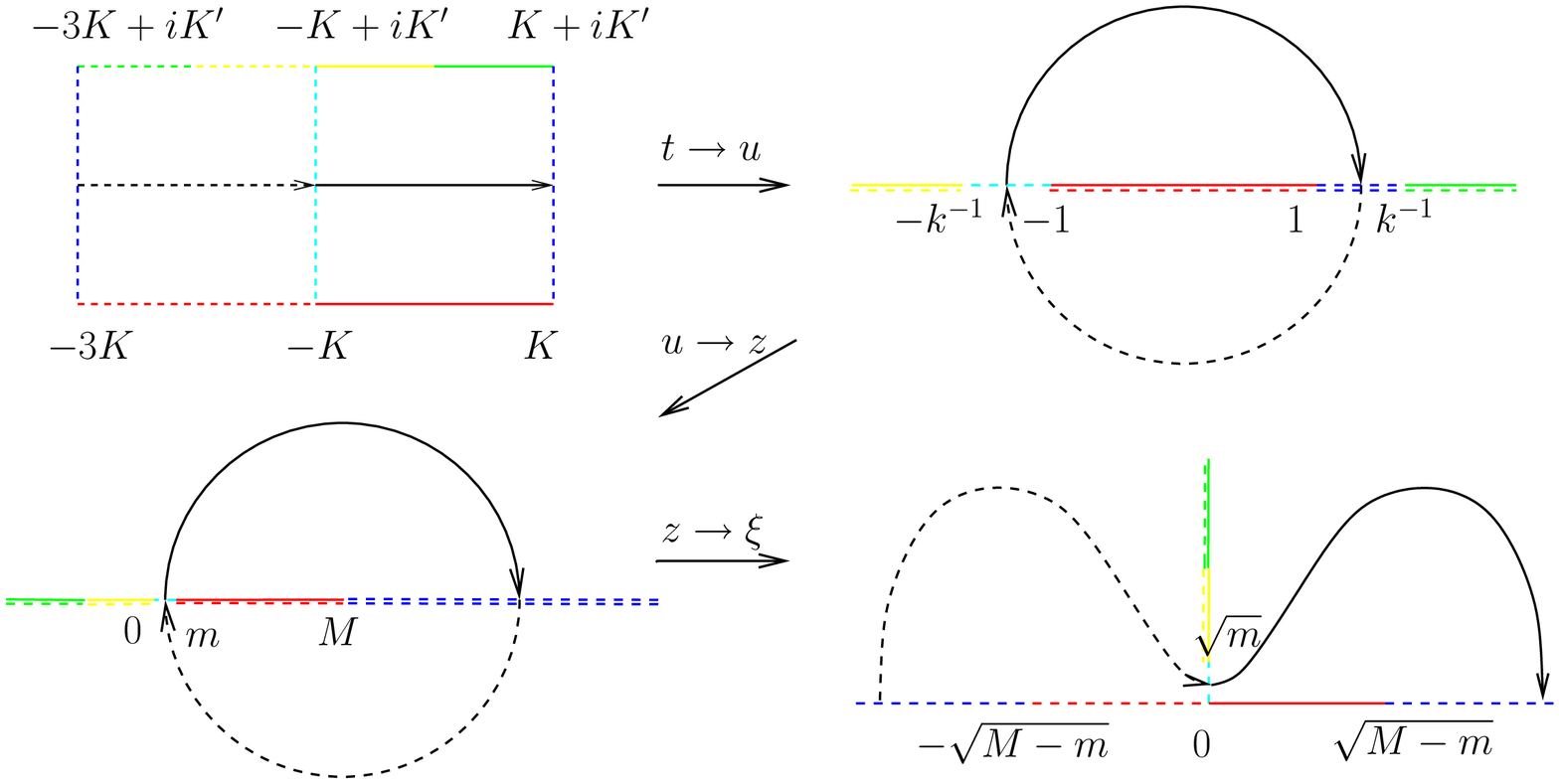}
  \end{center}
  \caption{The map from the rectangular domain $[-3K, K]\times[0, K']$
    to the upper-half of the domain $U$. The map is constructed in
    three steps: $t \to u \to z \to \xi$. The boundaries are shown in
    various colors and line styles.}
  \label{fig:map}
\end{figure}

Actually, what is done could be understood as follows.  Similar
to \cite{HaleHighamTrefethen2008}, we have constructed a map from the
rectangular domain $[-3K, K]\times[0, K']$ to the upper half of the
domain
\begin{equation*}
  U = \{z \mid \Im z \geq 0\} \backslash \bigl([-E_M, E_M]\cup 
  i[\pi/\beta, \infty)\bigr).
\end{equation*}
The map is carried out in three steps, shown in Fig.~\ref{fig:map}.
The first two steps use the original map constructed in
\cite{HaleHighamTrefethen2008}, however with extended domain of
definition. First, the Jacobi elliptic function
\begin{equation}
u = \sn(t) = \sn(t|k),\quad k = \frac{\sqrt{M/m} - 1}{\sqrt{M/m} + 1}
\end{equation}
maps the rectangular domain to the complex plane, with the ends
mapping to $[1,k^{-1}]$ and the middle vertical line $-K + i[0, K']$
to $[-k^{-1}, -1]$. Then, the M\"obius transformation
\begin{equation}
z = \sqrt{mM}\biggl(\frac{k^{-1}+u}{k^{-1}-u} \biggr)
\end{equation}
maps the complex plane to itself in such a way that
$[-k^{-1}, -1]$ and $[1,k^{-1}]$ are mapped to $[0,m]$ and
$[M,\infty]$, respectively. Finally, the shifted square root function
\begin{equation}
\xi = (z-m)^{1/2}
\end{equation}
maps the complex plane to the upper-half plane (we choose the branch
of the square root such that the lower-half plane is mapped to the
second quadrant and the upper-half plane is mapped to the first
quadrant), in such a way that $[0,m]$ is sent to $i[0,
\sqrt{m}]$ and $[M,\infty)$ is sent to $(-\infty, -\sqrt{M-m}] \cup
[\sqrt{M-m},\infty)$. The map can be extended to a map  from $[-7K, K]\times
[0,K']$ to the whole $U$, in this case, the $z$-plane becomes a
double-covered Riemann surface with branch point at $m$. 

Since the function $g$ is analytic in the domain $U$, the composite
function $g(t) = g(\xi(z(u(t))))$ is analytic in the stripe in the
$t$-plane, and therefore, the trapezoidal rule converges exponentially 
fast.
Using a similar analysis that leads to \eqref{eq:accuracy}, it can be
shown that
\begin{equation}
\label{eqn:gaplesserror}
\lVert g(\boldsymbol{B}) - g_Q(\boldsymbol{B})\rVert
=\Or(e^{- C Q/\log(\beta E_M)}),
\end{equation}
where $C$ is a constant.

We remark that the construction proposed in this subsection also
applies to the gapped case. In practice, if the temperature is high
(so that $\beta$ is small) or the gap around the chemical potential is
small (in particular, for gapless system), the contour passing through
the imaginary axis will be favorable; otherwise, the construction in
the last subsection will be more efficient.

\section{Rational approximations based on the multipole expansion}
\label{sec:FMM}

Another strategy for obtaining an efficient rational approximation for the
Fermi-Dirac function for finite temperature is based on the multipole
expansion, proposed recently in \cite{LinLuCarE2009}.  Let us first
recall the construction of the multipole representation.

Using the Matsubara representation (pole expansion) of the Fermi-Dirac
function, the density matrix is given by
\begin{equation}\label{eq:FOE}
\bvec{\rho} = 1 - 4\Re \sum_{l=1}^{\infty} \frac{1}{\beta(\bvec{H}-\mu) 
-(2l-1)\pi i}.
\end{equation}
The summation in \eqref{eq:FOE} can be seen as a summation of residues
contributed from the poles $\{(2l-1)\pi i\}$, with $l$ a positive
integer, on the imaginary axis. This suggests looking for a multipole
expansion of the contributions from the poles, as was done in the fast
multipole method (FMM) \cite{GreengardRokhlin1987}.  To do so, we use
a dyadic grouping of the poles, in which the $n$-th group
contains terms from $l=2^{n-1}$ to $l=2^{n}-1$, for a total of
$2^{n-1}$ terms. We decompose the summation in \eqref{eq:FOE}
accordingly. Let $x = \beta(\bvec{H}-\mu)$.  Then
\begin{equation}
\sum_{l=1}^{\infty} \frac{1}{x-(2l-1)\pi i}
=\sum_{n=1}^{\infty}\sum_{l=2^{n-1}}^{2^{n}-1} \frac{1}{x-(2l-1)\pi i}
=\sum_{n=1}^{\infty} S_n.
\end{equation}
The basic idea is to combine the simple poles into a set of multipoles at $l =
l_n$, where $l_n$ is taken as the midpoint of the interval $[2^{n-1},2^n-1]$
\begin{equation}
l_n = \frac{3\cdot 2^{n-1}-1}{2}.
\end{equation}
Then the $S_n$ term in the above equation can be written as
\begin{equation}\label{eq:Sn}
\begin{aligned}
S_n &= \sum_{l=2^{n-1}}^{2^{n}-1}\frac{1}{x-(2l_n-1)\pi i - 2(l-l_n)\pi i}\\
&= \sum_{l=2^{n-1}}^{2^{n}-1}\frac{1}{x-(2l_n-1)\pi i} \sum_{\nu=0}^{\infty}
\Bigl(\frac{2(l-l_n)\pi i}{x-(2l_n-1)\pi i}\Bigr)^{\nu} \\
&= \sum_{l=2^{n-1}}^{2^{n}-1}\frac{1}{x-(2l_n-1)\pi i} \sum_{\nu=0}^{P-1}
\Bigl(\frac{2(l-l_n)\pi i}{x-(2l_n-1)\pi i}\Bigr)^{\nu} \\
&\qquad\qquad +\sum_{l=2^{n-1}}^{2^{n}-1}
\frac{1}{x-(2l-1)\pi i}\Bigl(\frac{2(l-l_n)\pi i}{x-(2l_n-1)\pi i}
\Bigr)^P.
\end{aligned}
\end{equation}
Using the fact that $x$ is real, the second term in
 \eqref{eq:Sn} can be bounded by
\begin{equation*}
\sum_{l=2^{n-1}}^{2^{n}-1} \Bigl\lvert\frac{1}{x-(2l-1)\pi i}\Bigr\rvert
\Bigl\lvert \frac{2(l-l_n)\pi i}{x-(2l_n-1)\pi i}
\Bigr\rvert^P 
\leq \sum_{l=2^{n-1}}^{2^{n}-1} \frac{1}{\lvert (2l-1)\pi\rvert}
\Bigl\lvert\frac{2(l-l_n)}{2l_n-1}\Bigr\rvert^P \leq \frac{1}{2\pi} 
\frac{1}{3^P}.
\end{equation*}
Therefore, if we approximate the sum $S_n$ by the first $P$ terms,
the error decays exponentially fast with $P$:
\begin{equation}
\label{eq:errorP}
\left\lvert S_n(x) - 
\sum_{l=2^{n-1}}^{2^{n}-1}\frac{1}{x-(2l_n-1)\pi i} \sum_{\nu=0}^{P-1}
\Bigl(\frac{2(l-l_n)\pi i}{x-(2l_n-1)\pi i}\Bigr)^{\nu}\right\rvert \leq 
\frac{1}{2\pi}\frac{1}{3^P}.
\end{equation}
The above analysis is of course standard from the view point of the
fast multipole method \cite{GreengardRokhlin1987}.  The overall
philosophy is also similar: given a preset error tolerance, one
selects the value of $P$, the number of terms to retain in $S_n$,
according to \eqref{eq:errorP}.

Moreover, the remainder of the sum in \eqref{eq:FOE} from
$l=M_{\pole}+1$ to $\infty$ has an explicit expression
\begin{equation}\label{eq:digamma}
\Re\sum_{l=M_{\pole}+1}^{\infty}\dfrac{1}{2x-(2l-1)i\pi} = \frac{1}{2\pi}\Im 
\psi\left(M_{\pole}+\frac{1}{2}+\frac{i}{\pi}x\right),
\end{equation}
where $\psi$ is the digamma function 
$\psi(z) = \Gamma'(z)/\Gamma(z)$. 

In summary, we arrive at the following multipole representation for
the Fermi operator \cite{LinLuCarE2009}:
\begin{multline}\label{eq:polesumfmm}
  \bvec{\rho} = 1 - 4 \Re\sum_{n=1}^{N_G}
  \sum_{l=2^{n-1}}^{2^{n}-1}\frac{1}{\beta(\bvec{H}-\mu)-(2l_n-1)\pi
    i} \sum_{\nu=0}^{P-1} \Bigl(\frac{2(l-l_n)\pi
    i}{\beta(\bvec{H}-\mu)-(2l_n-1)\pi i}\Bigr)^{\nu} \\
  -\frac{2}{\pi}\Im\psi
  \left(M_{\pole}+\frac{1}{2}+\frac{i}{2\pi}\beta(\bvec{H}-\mu)\right)
  +\Or(N_G/3^{P}).
\end{multline}
Here $N_G$ is the number of groups in the multipole representation.
$M_{\pole}=2^{N_{G}}-1$ is the number of poles that are effectively
represented in the original Matsubara representation. In practice,
$N_G$ simple poles are first calculated, and then the $N_G(P-1)$
multipoles can be constructed through matrix-matrix multiplication.

A disadvantage of \eqref{eq:polesumfmm} is that one needs to
multiply simple poles together to get the multipoles before extracting
the diagonal of Fermi operator. This prevents us from being able to directly
apply the fast algorithms for extracting the diagonal of an inverse
matrix, such as the one proposed in
\cite{LinLuYingCarE2009}. Therefore, it is useful to find an
expansion similar to \eqref{eq:polesumfmm} that uses only simple
poles. As we mentioned earlier, the key idea in deriving
\eqref{eq:polesumfmm} is to combine the poles in each group together
to form multipoles as the distance between them and the real axis is
large. However, if instead we want an expansion that involves only
simple poles, it is natural to revisit the variants of FMM that only
use simple poles, for example, the version introduced in
\cite{YingBirosZorin2004}. The basic idea there is to use a set of
equivalent charges on a circle surrounding the poles in each group to
reproduce the effective potential away from these poles.

Specifically, take the group of poles from $l = 2^{n-1}$ to $l = 2^n -
1$ for example. Consider a circle $B_n$ with center $c_n = (3\cdot
2^{n-1} - 2)\pi i$ and radius $r_n = 2^{n-1}\pi$. It is clear that the
circle $B_n$ encloses the poles considered. Take $P$ equally spaced
points $\{x_{n,k}\}_{k=1}^P$ on the circle $B_n$. Next, one needs to
place equivalent charges $\{\rho_{n,k}\}_{k=1}^P$ at these points such
that the potential produced by these equivalent charges match with the
potential produced by the poles inside $B_n$ away from the circle.
This can be done in several ways, for example, by matching the
multipole expansion, by discretizing the potential on $B_n$ generated
by the poles, and so on. Here we follow the approach used in
\cite{YingBirosZorin2004}. 

We simply take a bigger concentric circle $\BB_n$ outside $B_n$ with
radius $R_n = 2^n \pi$ and match the potential generated on $\BB_n$ by
the poles and by the equivalent charges on $B_n$.  For this purpose,
we solve for $\rho_{n,k}$ the equations
\begin{equation}
  \sum_{k=1}^P \frac{\rho_{n,k}}{y - x_{n,k}} 
  = \sum_{l=2^{n-1}}^{2^n -1} \frac{1}{y - (2l-1)\pi i}, \quad y \in \BB_n.
\end{equation}
Regularization techniques such as Tikhonov regularization are required
here since this is a first-kind Fredholm equation.

One can also prove that similar to the original version of the multipole
representation, the error in the potential produced by the equivalent
charges decay exponentially in $P$, the details can be found in
\cite{YingBirosZorin2004}. Putting these all together, we can write down
the following expansion of the Fermi-Dirac function
\begin{multline}\label{eq:polesumfmm2}
  \bvec{\rho} = 1 - 4 \Re\sum_{n=1}^{N_G}
  \sum_{k=1}^{P}\frac{\rho_{n,k}}{\beta(\bvec{H}-\mu)-x_{n,k}}\\
  -\frac{2}{\pi}\Im\psi
  \left(M_{\pole}+\frac{1}{2}+\frac{i}{2\pi}\beta(\bvec{H}-\mu)\right)
  +\Or(N_G/3^{P}).
\end{multline}
The number of poles that are effectively
represented in the original Matsubara representation is still
$M_{\pole}=2^{N_G}-1$. $N_{\pole}=N_{G}P$ simple poles are now to be
calculated in practice.

The tail part can be approximated using a Chebyshev polynomial
expansion.  Similar to the analysis in \cite{LinLuCarE2009}, it can be
shown that the complexity of the expansion is $\Or(\log\beta\Delta
E)$. As we pointed out earlier, the advantage of
\eqref{eq:polesumfmm2} over \eqref{eq:polesumfmm} is that only simple
poles are involved in the formula. This is useful when combined with
fast algorithms for extracting the diagonal of an inverse matrix
\cite{LinLuYingCarE2009}.

\begin{figure}[ht]
  \begin{center}
    \includegraphics[width=8cm]{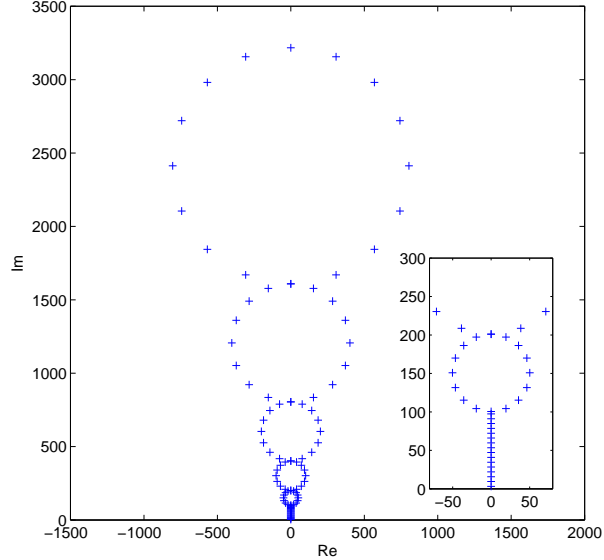}
  \end{center}
  \caption{A typical configuration of the poles in the multipole
    representation type algorithm.  $M_{\pole}=512$ and $P=16$ is used
    in this figure.  The poles with negative imaginary parts are not
    explicitly shown.  The inset shows the first few poles. The first
    $16$ poles are calculated separately and the starting level is
    $n=5$.}
  \label{fig:fmmcontour}
\end{figure}

Note that in \eqref{eq:polesumfmm} and \eqref{eq:polesumfmm2},
for $2^{n-1}<P$ there would be no savings if we use $P$ terms in the
expansion. They are written in this form just for simplicity. In
practice the first $P$ simple poles will be calculated separately and
the multipole expansion will be used starting from the $(P+1)$-th term
and the starting level is $n=\log_{2}P+1$. We show in
Fig.~\ref{fig:fmmcontour} a typical configuration of the set of poles
in the multipole representation type algorithm.

\section{Numerical results}
\label{sec:results}

We test the algorithms described above using a two dimensional nearest
neighbor tight binding model for the Hamiltonian. The matrix
components of the Hamiltonian can be written as (in atomic units),
\begin{equation}
  H_{i'j';ij} =  
  \begin{cases}
    2 + V_{ij}, &i'=i,j'=j,\\
    -1/2 + V_{ij}, &i'=i\pm 1, j'=j \text{ or } i'=i,j'=j\pm 1. 
  \end{cases}
  \label{eqn:Hamiltonian}
\end{equation}
The on-site potential energy $V_{ij}$ is chosen to be a uniform random
number between $0$ and $10^{-3}$. The domain size is $32\times 32$
with periodic boundary condition. The chemical potential will be 
specified later. 
The accuracy is measured by the $L^1$ error of the electronic density profile
per electron
\begin{equation}
  \Delta\rho_{\mathrm{rel}}=\frac{\tr{\abs{\widehat{P}-
  P} } }{N_{\mathrm{Electron}}}.
  \label{}
\end{equation}

\subsection{Contour integral representation: gapped case}

The error of the contour integral representation is determined by
$N_{\pole}$. At finite temperature $N_{\pole}=2Q$, while at zero
temperature $N_{\pole}=Q$, with $Q$ being the quadrature points on one
loop of the contour. The performance of the algorithm is studied by
the minimum number of $N_{\pole}$ such that
$\Delta\rho_{\mathrm{rel}}$ (the $L^1$ error in the electronic density
per electron) is smaller than $10^{-6}$.  For a given temperature, the
chemical potential $\mu$ is set to satisfy
\begin{equation}
  \Tr P = N_{\mathrm{Electron}}.
  \label{}
\end{equation}
In our setup the energy gap $E_{g}\approx 0.01\ \text{Hartree}=0.27\
\eV$ and $E_{M}\approx 4\ \text{Hartree}$. Therefore, this system can
be regarded as a crude model for semiconductor with a small energy
gap. The number of $N_{\pole}$ and the error
$\Delta\rho_{\mathrm{rel}}$ are shown in Table~\ref{tab:contourgap}
with respect to $\beta\Delta E$ ranging between $4,000$ and up to
$270,000$. Because of the existence of the finite energy gap, the
performance is essentially independent of $\beta\Delta E$, as is
clearly shown in Table~\ref{tab:contourgap}.

\begin{table}
  \centering
  \begin{tabular}{c||c|c}
    \toprule
    $\beta\Delta E$ & $N_{\pole}$ & 
    $\Delta\rho_{\mathrm{rel}}$ \\
    \midrule
    $4,208$ & $40$ & $5.68\times 10^{-7}$ \\
    $8,416$ & $44$ & $3.86\times 10^{-7}$ \\
    $16,832$ & $44$ & $3.60\times 10^{-7}$ \\
    $33,664$ & $44$ & $3.55\times 10^{-7}$ \\
    $67,328$ & $44$ & $3.57\times 10^{-7}$ \\
    $134,656$ & $44$ & $3.47\times 10^{-7}$ \\
    $269,312$ & $44$ & $3.55\times 10^{-7}$ \\
    \bottomrule
  \end{tabular}
  \vspace{1em}
  \caption{$N_{\pole}$ and $L^{1}$ error of electronic density per electron
    with respect to various $\beta\Delta E$. The energy gap $E_{g}\approx 0.01$.
    The contour integral representation for gapped system at finite temperature
    is  used for the calculation.  
    The performance of the algorithm depends weakly on $\beta\Delta E$.}
  \label{tab:contourgap}
\end{table}

\FloatBarrier

When the temperature is low and therefore when $\beta$ is large, as
discussed before the finite temperature result is well approximated by the
zero temperature Fermi operator, \ie, the matrix sign function. In such case
the quadrature formula is given by \eqref{eqn:matrixsignquad}. Only the
contour that encircles the spectrum lower than chemical potential is
calculated, and $N_{\pole} = Q$.  

\begin{figure}[ht]
  \begin{center}
    \includegraphics[width=8cm]{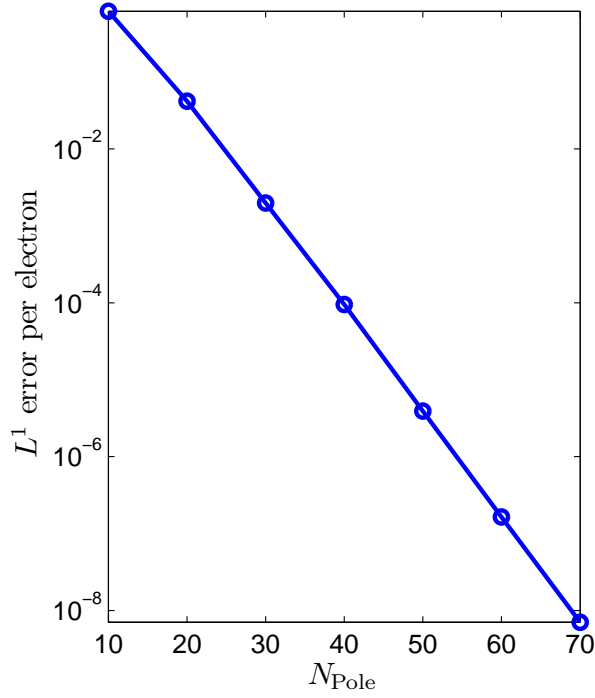}
  \end{center}
  \caption{The lin-log plot of the $L^{1}$ error of electronic density
    per electron with respect to $N_{\pole}$. The energy gap
    $E_{g}\approx 10^{-6}$. The contour integral representation for
    gapped system at zero-temperature is used for calculation. }
  \label{fig:matrixsign}
\end{figure}

In order to study the dependence of $\Delta\rho_{\mathrm{rel}}$ on the
number of poles $N_{\pole}$, we tune artificially the chemical
potential to reduce the energy gap to $10^{-6}\ \text{Hartree}$.
Fig.~\ref{fig:matrixsign} shows the exponential decay of
$\Delta\rho_{\mathrm{rel}}$ with respect to $N_{\pole}$.  For example,
in order to reach the $10^{-6}$ error criterion, $N_{\pole}\approx 50$
is sufficient.  The increase in $N_{\pole}$ is very small compared to the
large decrease of energy gap and this is consistent the logarithmic
dependence of $N_{\pole}$ on $E_g$ given by \eqref{eq:gapaccuracy}.

\FloatBarrier

\subsection{Contour integral representation: gapless case}

\begin{figure}[ht]
  \begin{center}
    \includegraphics[width=8cm]{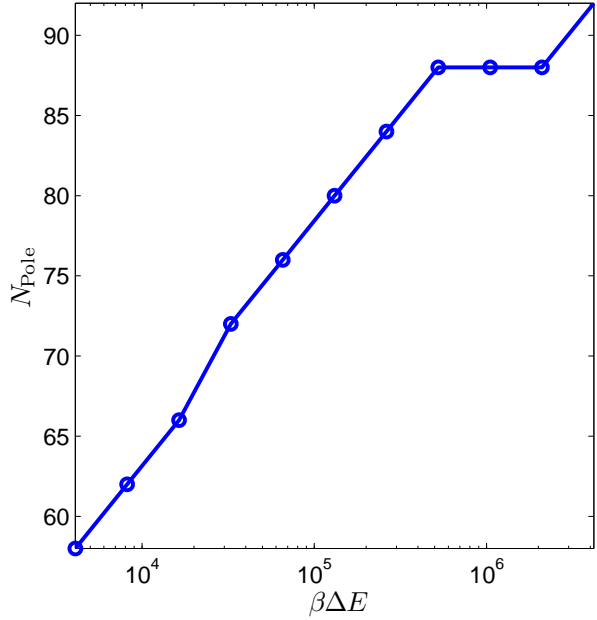}
  \end{center}
  \caption{Log-lin plot of $N_{\pole}$ with respect to $\beta\Delta
    E$.  The contour integral representation for gapless system is
    used for the calculation. }
  \label{fig:contournbetae}
\end{figure}

For gapless systems such as metallic systems, our quadrature
formula in \eqref{eqn:gapless} exploits the effective gap on the
imaginary axis due to finite temperature. In the following
results the chemical potential is set artificially so that
$E_{g}=0$. $E_{M}\approx 4\ \text{Hatree}$ and the error criterion is
still $10^{-6}$ as in the gapped case. Table~\ref{tab:contourgapless}
reports the number of poles $N_{\pole}$ and the error
$\Delta\rho_{\mathrm{rel}}$ with respect to $\beta\Delta E$ ranging
from $4,000$ up to $4$ million. These results are further summarized
in Fig.~\ref{fig:contournbetae} to show the logarithmic dependence of
$N_{\pole}$ on $\beta\Delta E$, as predicted in the analysis of
\eqref{eqn:gaplesserror}.

\begin{table}
  \centering
  \begin{tabular}{c||c|c}
    \toprule
    $\beta\Delta E$ & $N_{\pole}$ & 
    $\Delta\rho_{\mathrm{rel}}$ \\
    \midrule
    $4,208$ & $58$ & $1.90\times 10^{-7}$ \\
    $8,416$ & $62$ & $5.32\times 10^{-7}$ \\
    $16,832$ & $66$ & $8.28\times 10^{-7}$ \\
    $33,664$ & $72$ & $3.55\times 10^{-7}$ \\
    $67,328$ & $76$ & $3.46\times 10^{-7}$ \\
    $134,656$ & $80$ & $1.69\times 10^{-7}$ \\
    $269,312$ & $84$ & $8.89\times 10^{-8}$ \\
    $538,624$ & $88$ & $7.09\times 10^{-8}$ \\
    $1,077,248$ & $88$ & $8.94\times 10^{-7}$ \\
    $2,154,496$ & $88$ & $4.25\times 10^{-7}$ \\
    $4,308,992$ & $92$ & $3.43\times 10^{-7}$ \\
    \bottomrule
  \end{tabular}
  \vspace{1em}
  \caption{$N_{\pole}$ and $L^{1}$ error of electronic density per electron
    with respect to various $\beta\Delta E$. $E_{g}=0$. The contour integral
    representation for gapless system is used for the
    calculation.
  } \label{tab:contourgapless}
\end{table}

\FloatBarrier

\subsection{Multipole representation}

The approach \eqref{eq:polesumfmm2} based on the multipole
representation has three parts of error: the finite-term multipole
expansion, the finite-term Chebyshev expansion for the tail part, and
the truncated matrix-matrix multiplication in the Chebyshev expansion.

The error from the multipole expansion is well controlled by $P$ in
\eqref{eq:polesumfmm2}. When $P=16$, $1/3^P\sim \Or(10^{-8})$. The
number of groups $N_G$ is usually no more than $20$, and therefore the
error introduced by multipole expansion is around $\Or(10^{-7})$,
which is much less than the error criterion $10^{-6}$.

The number of terms in the Chebyshev expansion for the tail part
$N_{\mathrm{Cheb}}$ is $\Or\bigl( \frac{\beta\Delta E}{M_{\pole}}
\bigr)$, with $M_{\pole}$ being the number of poles that are excluded
in the tail part in the pole expansion. The truncation radius for the
tail part is $\Or\bigl( \exp(-C\frac{\beta\Delta
  E}{M_{\pole}})\bigr)$. In order to reach a fixed target accuracy, we
set $M_{\pole}$ to be proportional to $\beta\Delta E$. Due to the fact
that $M_{\pole}\approx 2^{N_G} \approx 2^{N_{\pole}/P}$, this requires
$N_{\pole}$ to grow logarithmically with respect to $\beta\Delta E$.

\begin{figure}[ht]
  \begin{center}
    \includegraphics[width=8cm]{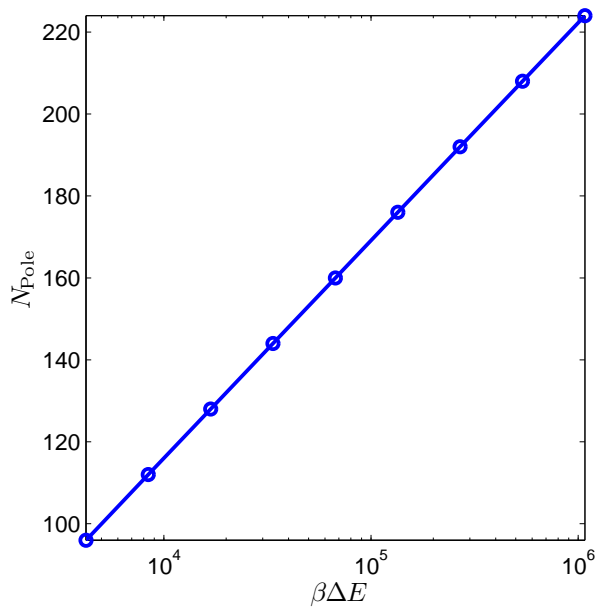}
  \end{center}
	\caption{log-lin plot of $N_{\pole}$ with respect to $\beta\Delta E$. The
	multipole representation is used for the calculation. }
  \label{fig:fmmnbetae}
\end{figure}

The target accuracy for the Chebyshev expansion is set to be
$10^{-7}$ and the truncation radius for the tail is set to be $4$ for
the metallic system under consideration. For $\beta\Delta E=4208$,
$M_{\pole}$ is set to be $512$ so that the error is smaller than
$10^{-6}$. For other cases, $M_{\pole}$ scales linearly with
$\beta\Delta E$. The lin-log plot in Fig.~\ref{fig:fmmnbetae} shows
the logarithmic dependence of $N_{\pole}$ with respect to $\beta\Delta
E$.  For more detailed results, Table~\ref{tab:multipole} measures
$M_{\pole}$, $N_{\pole}$, $N_{\mathrm{Cheb}}$, and
$\Delta\rho_{\mathrm{rel}}$ for $\beta\Delta E$ ranging from $4000$ up
to $1$ million.  For all cases, $N_{\mathrm{Cheb}}$ is kept as a small
constant. Note that the truncation radius is always set to be a small
number $4$, and this indicates the tail part is extremely localized in
the multipole representation due to the effectively raised
temperature.

Table~\ref{tab:multipole} indicates that the error exhibits some
slight growth. We believe that it comes from the growth of the number
of groups in the multipole representation \eqref{eq:polesumfmm2} and
also the extra $\log\log$ dependence on $\beta\Delta E$ (see
\cite{LinLuCarE2009} for details). When compared with the results
reported in Table~\ref{tab:contourgapless}, we see that for the
current application to electronic structure, the contour integral
representation outperforms the multipole representation in terms of
both the accuracy and the number of poles used.

\begin{table}
  \centering
  \begin{tabular}{c||c|c|c|c}
    \toprule
    $\beta\Delta E$ & $M_{\pole}$ & $N_{\pole}$ & $N_{\mathrm{Cheb}}$ &
    $\Delta\rho_{\mathrm{rel}}$ \\
    \midrule
    $4,208$ & $512$ & $96$ & $22$ & $4.61\times 10^{-7}$ \\
    $8,416$ & $1,024$ & $112$ & $22$ & $4.76\times 10^{-7}$ \\
    $16,832$ & $2,048$ & $128$ & $22$ & $4.84\times 10^{-7}$ \\
    $33,664$ & $4,096$ & $144$ & $22$ & $4.88\times 10^{-7}$ \\
    $67,328$ & $8,192$ & $160$ & $22$ & $4.90\times 10^{-7}$ \\
    $134,656$ & $16,384$ & $176$ & $22$ & $4.90\times 10^{-7}$ \\
    $269,312$ & $32,768$ & $192$ & $22$ & $6.98\times 10^{-7}$ \\
    $538,624$ & $65,536$ & $208$ & $22$ & $3.20\times 10^{-6}$ \\
    $1,077,248$ & $131,072$ & $224$ & $22$ & $7.60\times 10^{-6}$  \\
    \bottomrule
  \end{tabular}
  \vspace{1em}
  \caption{The number of poles calculated $N_{\pole}$, the order of Chebyshev
	expansion for the tail part $N_{\mathrm{Cheb}}$, and the $L^1$ error of
	electronic density per electron with respect to various $\beta\Delta E$.
	The number of poles excluded in the tail part $M_{\pole}$ is chosen to
	be proportional to $\beta\Delta E$.}
  \label{tab:multipole}
\end{table}

\section{Conclusion}

We propose two approaches for the expansion of Fermi operator: a
rational approximation based on the contour integral idea introduced
in \cite{HaleHighamTrefethen2008} and a variant of the multipole
representation in \cite{LinLuCarE2009} using only simple poles. Both
approximations result in logarithmic scaling complexity with respect
to $\beta\Delta\eps$ with small prefactor. Fast algorithms for
electronic structure calculations can be obtained by combining these
approaches with the algorithm introduced in \cite{LinLuYingCarE2009}
for extracting the diagonal of the inverse of a matrix.

\vspace{1em}
\noindent{\bf Acknowledgement:}
This is a continuation of the work that was done jointly with 
Roberto Car, to whom we are very grateful for many stimulating discussions.
This work was partially supported by DOE under Contract
No. DE-FG02-03ER25587 and by ONR under Contract No. N00014-01-1-0674
(L.~L., J.~L. and W.~E), and by an Alfred P. Sloan Research Fellowship
and a startup grant from University of Texas at Austin (L.~Y.).

\bibliographystyle{abbrv}
\bibliography{contour}

\begin{thebibliography}{10}

\bibitem{BaroniGiannozzi1992}
S.~Baroni and P.~Giannozzi.
\newblock Towards very large-scale electronic-structure calculations.
\newblock {\em Europhys. Lett.}, 17(6):547--552, 1992.

\bibitem{parrinelloarxiv}
M.~Ceriotti, T.~K{\"u}hne, and M.~Parrinello.
\newblock {A hybrid approach to Fermi operator expansion}.
\newblock {\em arXiv:0809.2232v1}, 2008.

\bibitem{parrinello2008}
M.~Ceriotti, T.~K{\"u}hne, and M.~Parrinello.
\newblock {An efficient and accurate decomposition of the Fermi operator.}
\newblock {\em J. Chem. Phys}, 129(2):024707, 2008.

\bibitem{goedecker1999}
S.~Goedecker.
\newblock {Linear scaling electronic structure methods}.
\newblock {\em Rev. Mod. Phys.}, 71(4):1085--1123, 1999.

\bibitem{goedeckercolombo1994}
S.~Goedecker and L.~Colombo.
\newblock Efficient linear scaling algorithm for tight-binding molecular
  dynamics.
\newblock {\em Phys. Rev. Lett.}, 73(1):122--125, Jul 1994.

\bibitem{goedeckerteter1995}
S.~Goedecker and M.~Teter.
\newblock Tight-binding electronic-structure calculations and tight-binding
  molecular dynamics with localized orbitals.
\newblock {\em Phys. Rev. B}, 51(15):9455--9464, Apr 1995.

\bibitem{GreengardRokhlin1987}
L.~Greengard and V.~Rokhlin.
\newblock A fast algorithm for particle simulations.
\newblock {\em J. Comput. Phys.}, 73(2):325--348, 1987.

\bibitem{HaleHighamTrefethen2008}
N.~Hale, N.~J. Higham, and L.~N. Trefethen.
\newblock Computing \protect{$A^{\alpha}$}, \protect{$\log(A)$}, and related
  matrix functions by contour integrals.
\newblock {\em SIAM J. Numer. Anal.}, 46(5):2505--2523, 2008.

\bibitem{parrinello2005}
F.~Krajewski and M.~Parrinello.
\newblock Stochastic linear scaling for metals and nonmetals.
\newblock {\em Phys. Rev. B}, 71(23):233105, 2005.

\bibitem{headgordon2004}
W.~Liang, R.~Baer, C.~Saravanan, Y.~Shao, A.~T. Bell, and M.~Head-Gordon.
\newblock Fast methods for resumming matrix polynomials and {C}hebyshev matrix
  polynomials.
\newblock {\em J. Comput. Phys.}, 194(2):575 -- 587, 2004.

\bibitem{headgordon2003}
W.~Liang, C.~Saravanan, Y.~Shao, R.~Baer, A.~T. Bell, and M.~Head-Gordon.
\newblock Improved {F}ermi operator expansion methods for fast electronic
  structure calculations.
\newblock {\em J. Chem. Phys.}, 119(8):4117--4125, 2003.

\bibitem{LinLuCarE2009}
L.~Lin, J.~Lu, R.~Car, and W.~E.
\newblock Multipole representation of the {F}ermi operator with application to
  the electronic structure analysis of metallic systems.
\newblock {\em Phys. Rev. B}, 79:115133, 2009.

\bibitem{LinLuYingCarE2009}
L.~Lin, J.~Lu, L.~Ying, R.~Car, and W.~E.
\newblock Fast algorithm for extracting the diagonal of the inverse matrix with
  application to the electronic structure analysis of metallic systems.
\newblock {\em submitted}, 2009.

\bibitem{Ozaki2007}
T.~Ozaki.
\newblock Continued fraction representation of the {F}ermi-{D}irac function for
  large-scale electronic structure calculations.
\newblock {\em Phys. Rev. B}, 75(3):035123, 2007.

\bibitem{YingBirosZorin2004}
L.~Ying, G.~Biros, and D.~Zorin.
\newblock A kernel-independent adaptive fast multipole algorithm in two and
  three dimensions.
\newblock {\em J. Comput. Phys.}, 196(2):591 -- 626, 2004.

\end{thebibliography}

\end{document}